# Approximation to real numbers by cubic algebraic integers II

by Damien ROY

**Abstract.** It has been conjectured for some time that, for any integer $n \geq 2$, any real number $\epsilon > 0$ and any transcendental real number $\xi$, there would exist infinitely many algebraic integers $\alpha$ of degree at most $n$ with the property that $|\xi - \alpha| \leq H(\alpha)^{-n+\epsilon}$, where $H(\alpha)$ denotes the height of $\alpha$. Although this is true for $n = 2$, we show here that, for $n = 3$, the optimal exponent of approximation is not 3 but $(3 + \sqrt{5})/2 \simeq 2.618$.

## 1. Introduction

Define the *height* $H(\alpha)$ of an algebraic number $\alpha$ as the largest absolute value of the coefficients of its irreducible polynomial over **Z**. Thanks to work of H. Davenport and W. M. Schmidt, we know that, for any real number $\xi$ which is neither rational nor quadratic over **Q**, there exists a constant $c > 0$ such that the inequality

$$|\xi - \alpha| \leq cH(\alpha)^{-\gamma^2},$$

where $\gamma = (1 + \sqrt{5})/2$ denotes the golden ratio, has infinitely many solutions in algebraic integers $\alpha$ of degree at most 3 over **Q** (see Theorem 1 of [3]). The purpose of this paper is to show that the exponent $\gamma^2$ in this statement is best possible.

**Theorem 1.1.** *There exists a real number $\xi$ which is transcendental over **Q** and a constant $c_1 > 0$ such that, for any algebraic integer $\alpha$ of degree at most 3 over **Q**, we have*

$$|\xi - \alpha| \geq c_1 H(\alpha)^{-\gamma^2}.$$

In general, for a positive integer $n$, denote by $\tau_n$ the supremum of all real numbers $\tau$ with the property that any transcendental real number $\xi$ admits infinitely many approximations by algebraic integers $\alpha$ of degree at most $n$ over **Q** with $|\xi - \alpha| \leq H(\alpha)^{-\tau}$. Then, the above result, shows that $\tau_3 = \gamma^2 \simeq 2.618$ against the natural conjecture that $\tau_n = n$ for all $n \geq 2$ (see [7], page 259). Since $\tau_2 = 2$ (see the introduction of [3]), it leaves open the problem of evaluating $\tau_n$ for $n \geq 4$. At present the best known estimates valid for general $n \geq 2$ are

$$\lceil (n+1)/2 \rceil \leq \tau_n \leq n$$

where the upper bound comes from standard metrical considerations, while the lower bound, due to M. Laurent [4], refines, for even integers $n$, the preceding lower bound

Work partly supported by NSERC and CICMA.

$\tau_n \geq \lfloor (n+1)/2 \rfloor$ of Davenport and Schmidt [3]. Note that similar estimates are known for the analog problem of approximation by algebraic numbers, but in this case the optimal exponent is known only for $n \leq 2$ (see [2]).

In the next section we recall the results that we will need from [6]. Then, in Section 3, we present the class of real numbers for which we will prove, in Section 4, that they satisfy the measure of approximation of Theorem 1.1. Section 3 also provides explicit examples of such numbers based on the Fibonacci continued fractions of [5] and [6] (a special case of the Sturmian continued fractions of [1]).

## 2. Extremal real numbers

The arguments of Davenport and Schmidt in Section 2 of [3] show that, if a real number $\xi$ is not algebraic over $\mathbf{Q}$ of degree at most 2 and has the property stated in Theorem 1.1, then there exists another constant $c_2 > 0$ such that the inequalities

$$1 \leq |x_0| \leq X, \quad |x_0\xi - x_1| \leq c_2 X^{-1/\gamma}, \quad |x_0\xi^2 - x_2| \leq c_2 X^{-1/\gamma}, \tag{2.1}$$

have a solution in integers $x_0, x_1, x_2$ for any real number $X \geq 1$. In [6], we defined a real number $\xi$ to be *extremal* if it is not algebraic over $\mathbf{Q}$ of degree at most 2 and satisfies the latter property of simultaneous approximation. We showed that such numbers exist and form an enumerable set. Thus, candidates for Theorem 1.1 have to be extremal real numbers.

For each $\mathbf{x} = (x_0, x_1, x_2) \in \mathbf{Z}^3$ and each $\xi \in \mathbf{R}$, we define

$$\|\mathbf{x}\| = \max\{|x_0|, |x_1|, |x_2|\} \quad \text{and} \quad L_\xi(\mathbf{x}) = \max\{|x_0\xi - x_1|, |x_0\xi^2 - x_2|\}.$$

Identifying $\mathbf{x}$ with the symmetric matrix

$$\begin{pmatrix} x_0 & x_1 \\ x_1 & x_2 \end{pmatrix},$$

we also define

$$\det(\mathbf{x}) = x_0 x_2 - x_1^2.$$

Then, Theorem 5.1 of [6] provides the following characterization of extremal real numbers.

**Proposition 2.1.** *A real number $\xi$ is extremal if and only if there exists a constant $c_3 \geq 1$ and an unbounded sequence of non-zero points $(\mathbf{x}_k)_{k \geq 1}$ of $\mathbf{Z}^3$ satisfying, for all $k \geq 1$,*
  *(i) $c_3^{-1}\|\mathbf{x}_k\|^\gamma \leq \|\mathbf{x}_{k+1}\| \leq c_3\|\mathbf{x}_k\|^\gamma$,*
  *(ii) $c_3^{-1}\|\mathbf{x}_k\|^{-1} \leq L_\xi(\mathbf{x}_k) \leq c_3\|\mathbf{x}_k\|^{-1}$,*
  *(iii) $1 \leq |\det(\mathbf{x}_k)| \leq c_3$,*
  *(iv) $1 \leq |\det(\mathbf{x}_k, \mathbf{x}_{k+1}, \mathbf{x}_{k+2})| \leq c_3$.*

In order to prove our main Theorem 1.1, we will also need the following special case of Proposition 9.1 of [6] where, for a real number $t$, the symbol $\{t\}$ denotes the distance from $t$ to a closest integer:

**Proposition 2.2.** *Let $\xi$ be an extremal real number and let $(\mathbf{x}_k)_{k\geq 1}$ be as in Proposition 2.1. Assume that, upon writing $\mathbf{x}_k = (x_{k,0}, x_{k,1}, x_{k,2})$, there exists a constant $c_4 > 0$ such that*

$$\{x_{k,0}\xi^3\} \geq c_4$$

*for all $k \geq 1$. Then, for any algebraic integer $\alpha$ of degree at most 3 over $\mathbf{Q}$, we have*

$$|\xi - \alpha| \geq c_5 H(\alpha)^{-\gamma^2}$$

*for some constant $c_5 > 0$.*

Since extremal real numbers are transcendental over $\mathbf{Q}$ (see §5 of [6]), this reduces the proof of Theorem 1.1 to finding extremal real numbers satisfying the hypotheses of the above proposition. Note that, for an extremal real number $\xi$ and a corresponding sequence $(\mathbf{x}_k)_{k\geq 1}$, Proposition 9.2 of [6] shows that there exists a constant $c_6 > 0$ such that

$$\{x_{k,0}\xi^3\} \geq c_6 \|\mathbf{x}_k\|^{-1/\gamma^3}$$

for any sufficiently large $k$.

We also mention the following direct consequence of Corollary 5.2 of [6]:

**Proposition 2.3.** *Let $\xi$ be an extremal real number and let $(\mathbf{x}_k)_{k\geq 1}$ be as in Proposition 2.1. Then there exists an integer $k_0 \geq 1$ and a $2 \times 2$ matrix $M$ with integral coefficients such that, viewing each $\mathbf{x}_k$ as a symmetric matrix, the point $\mathbf{x}_{k+2}$ is a rational multiple of $\mathbf{x}_{k+1}M\mathbf{x}_k$ when $k \geq k_0$ is odd, and a rational multiple of $\mathbf{x}_{k+1}{}^t M\mathbf{x}_k$ when $k \geq k_0$ is even.*

*Proof.* Corollary 5.2 together with formula (2.2) of [6] show that there exists an integer $k_0 \geq 1$ such that $\mathbf{x}_{k+2}$ is a rational multiple of $\mathbf{x}_{k+1}\mathbf{x}_{k-1}^{-1}\mathbf{x}_{k+1}$ for all $k > k_0$. If $S$ is a $2 \times 2$ matrix such that $\mathbf{x}_{k+1}$ is a rational multiple of $\mathbf{x}_k S\mathbf{x}_{k-1}$ for some $k > k_0$, this implies that $\mathbf{x}_{k+2}$ is a rational multiple of $\mathbf{x}_k S\mathbf{x}_{k+1}$ and thus, by taking transpose, that $\mathbf{x}_{k+2}$ is a rational multiple of $\mathbf{x}_{k+1}{}^t S\mathbf{x}_k$. The conclusion then follows by induction on $k$, upon choosing $M$ so that the required property holds for $k = k_0$. □

Note that, in the case where all points $\mathbf{x}_k$ have determinant 1, one may assume that $M \in \mathrm{GL}_2(\mathbf{Z})$ in the above proposition and the conclusion then becomes $\mathbf{x}_{k+1} = \pm\mathbf{x}_{k+1}S\mathbf{x}_k$ where $S$ is either $M$ or ${}^t M$ depending on the parity of $k \geq k_0$. This motivates the following definition:

**Definition 2.4.** *Let $M \in \mathrm{GL}_2(\mathbf{Z})$ be a non-symmetric matrix. We denote by $\mathcal{E}(M)$ the set of extremal real numbers $\xi$ with the following property. There exists a sequence of points $(\mathbf{x}_k)_{k\geq 1}$ in $\mathbf{Z}^3$ satisfying the conditions of Proposition 2.1 which, viewed as symmetric matrices, belong to $\mathrm{GL}_2(\mathbf{Z})$ and verify the recurrence relation*

$$\mathbf{x}_{k+2} = \mathbf{x}_{k+1}S\mathbf{x}_k, \quad (k \geq 1), \quad \text{where} \quad S = \begin{cases} M & \text{if } k \text{ is odd,} \\ {}^t M & \text{if } k \text{ is even.} \end{cases}$$

Examples of extremal real numbers are the Fibonacci continued fractions $\xi_{a,b}$ (see [5] and §6 of [6]) where $a$ and $b$ denote distinct positive integers. They are defined as the real numbers

$$\xi_{a,b} = [0, a, b, a, a, b, \ldots] = 1/(a + 1/(b + \cdots))$$

whose sequence of partial quotients begins with 0 followed by the elements of the Fibonacci word on $\{a, b\}$, the infinite word $abaab\cdots$ starting with $a$ which is a fixed point of the substitution $a \mapsto ab$ and $b \mapsto a$. Corollary 6.3 of [6] then shows that such a number $\xi_{a,b}$ belongs to $\mathcal{E}(M)$ with

$$(2.5) \qquad M = \begin{pmatrix} a & 1 \\ 1 & 0 \end{pmatrix} \begin{pmatrix} b & 1 \\ 1 & 0 \end{pmatrix} = \begin{pmatrix} ab+1 & a \\ b & 1 \end{pmatrix}.$$

We conclude this section with the following result.

**Lemma 2.6.** *Assume that $\xi$ belongs to $\mathcal{E}(M)$ for some non-symmetric matrix*

$$M = \begin{pmatrix} a & b \\ c & d \end{pmatrix} \in \mathrm{GL}_2(\mathbf{Z}),$$

*and let $(\mathbf{x}_k)_{k \geq 1}$ be as in Definition 2.4. Then, upon writing $\mathbf{x}_k = (x_{k,0}, x_{k,1}, x_{k,2})$, we have, for all $k \geq 2$,*

(i) $\mathbf{x}_{k+2} = (ax_{k,0} + (b+c)x_{k,1} + dx_{k,2})\mathbf{x}_{k+1} \pm \mathbf{x}_{k-1}$,

(ii) $x_{k,0}x_{k+1,2} - x_{k,2}x_{k+1,0} = \pm(ax_{k-1,0} - dx_{k-1,2}) \pm (b-c)x_{k-1,1}$

*Proof.* For $k \geq 1$, we have

$$\mathbf{x}_{k+1} = \mathbf{x}_k S \mathbf{x}_{k-1} \quad \text{and} \quad \mathbf{x}_{k+2} = \mathbf{x}_{k+1}{}^t S \mathbf{x}_k$$

where $S$ is $M$ or ${}^tM$ according to whether $k$ is even or odd, and so

$$\mathbf{x}_{k+2} = {}^t\mathbf{x}_{k+2} = \mathbf{x}_k S \mathbf{x}_{k+1} = (\mathbf{x}_k S)^2 \mathbf{x}_{k-1}.$$

Since Cayley-Hamilton's theorem gives

$$(\mathbf{x}_k S)^2 = \mathrm{trace}(\mathbf{x}_k S)\mathbf{x}_k S - \det(\mathbf{x}_k S)I,$$

we deduce

$$\mathbf{x}_{k+2} = \mathrm{trace}(\mathbf{x}_k S)\mathbf{x}_{k+1} - \det(\mathbf{x}_k S)\mathbf{x}_{k-1}$$

which proves (i). Finally, (ii) follows from the fact that the left hand side of this inequality is the sum of the coefficients outside of the diagonal of the product

$$\mathbf{x}_k J \mathbf{x}_{k+1} \quad \text{where} \quad J = \begin{pmatrix} 0 & 1 \\ -1 & 0 \end{pmatrix},$$

and that, since $J\mathbf{x}_k J = \pm \mathbf{x}_k^{-1}$, we have

$$\mathbf{x}_k J \mathbf{x}_{k+1} = \pm J\mathbf{x}_k^{-1}\mathbf{x}_{k+1} = \pm JS\mathbf{x}_{k-1}.$$

## 3. A smaller class of real numbers

Although we expect that all extremal real numbers $\xi$ satisfy a measure of approximation by algebraic integers of degree at most 3 which is close to that of Theorem 1.1, say with exponent $\gamma^2 + \epsilon$ for any $\epsilon > 0$, we could only prove in [6] that they satisfy a measure with exponent $(3/2)\gamma^2$ (see Theorem 1.5 of [6]). Here we observe that the formulas of Lemma 2.6 show a particularly simple arithmetic for the elements $\xi$ of $\mathcal{E}(M)$ when, in the notation of this lemma, the matrix $M$ has $b = 1$, $c = -1$ and $d = 0$. Taking advantage of this, we will prove:

**Theorem 3.1.** *Let $a$ be a positive integer. Then, any element $\xi$ of*

$$\mathcal{E}_a = \mathcal{E}\begin{pmatrix} a & 1 \\ -1 & 0 \end{pmatrix}$$

*satisfies the measure of approximation of Theorem 1.1.*

The proof of this result will be given in the next section. Below, we simply show that, for $a = 1$, the corresponding set of extremal real numbers is not empty.

**Proposition 3.2.** *Let $m$ be a positive integer. Then, the real number*

$$\eta = (m + 1 + \xi_{m,m+2})^{-1} = [0, m+1, m, m+2, m, m, m+2, \ldots]$$

*belongs to the set $\mathcal{E}_1$.*

*Proof.* We first note that, if a real number $\xi$ belongs to $\mathcal{E}(M)$ for some non-symmetric matrix $M \in \mathrm{GL}_2(\mathbf{Z})$ with corresponding sequence of symmetric matrices $(\mathbf{x}_k)_{k \geq 1}$, and if $C$ is any element of $\mathrm{GL}_2(\mathbf{Z})$, then the real number $\eta$ for which $(\eta, -1)$ is proportional to $(\xi, -1)C$ belongs to $\mathcal{E}(^tCMC)$ with corresponding sequence $(C^{-1}\mathbf{x}_k{}^tC^{-1})_{k \geq 1}$. The conclusion then follows since $\xi_{m,m+2}$ belongs to $\mathcal{E}(M)$ where $M$ is given by (2.5) with $a = m$ and $b = m+2$ and since

$$^tCMC = \begin{pmatrix} 1 & 1 \\ -1 & 0 \end{pmatrix} \quad \text{for} \quad C = \begin{pmatrix} 0 & -1 \\ -1 & m+1 \end{pmatrix}.$$

**Remark.** In fact, it can be shown that $\mathcal{E}_a$ is not empty for any integer $a \geq 1$. For example, consider the sequence of matrices $(\mathbf{x}_k)_{k \geq 1}$ defined recursively using the formula of Definition 2.4 with

$$\mathbf{x}_1 = \begin{pmatrix} 1 & 1 \\ 1 & 0 \end{pmatrix}, \quad \mathbf{x}_2 = \begin{pmatrix} a^3 + 2a & a^3 - a^2 + 2a - 1 \\ a^3 - a^2 + 2a - 1 & a^3 - 2a^2 + 3a - 2 \end{pmatrix} \quad \text{and} \quad M = \begin{pmatrix} a & 1 \\ -1 & 0 \end{pmatrix}.$$

Then, using similar arguments as in Section 6 of [6], it can be shown that $(\mathbf{x}_k)_{k\geq 1}$ is a sequence of symmetric matrices in $\mathrm{GL}_2(\mathbf{Z})$ which satisfies the four conditions of Proposition 2.2 for some real number $\xi$ which therefore belongs to $\mathcal{E}_a$.

## 4. Proof of Theorem 3.1

We fix a positive integer $a$, a real number $\xi \in \mathcal{E}_a$, and a corresponding sequence of points $(\mathbf{x}_k)_{k\geq 1}$ of $\mathbf{Z}^3$ as in Definition 2.4. For simplicity, we also define
$$X_k = \|\mathbf{x}_k\| \quad \text{and} \quad \delta_k = \{x_{k,2}\xi\}, \quad (k \geq 1).$$
The constant $c_3$ being as in Proposition 2.1, we first note that

(4.1)
$$\{x_{k,0}\xi\} \leq |x_{k,0}\xi - x_{k,1}| \leq c_3 X_k^{-1},$$
$$\{x_{k,1}\xi\} \leq |x_{k,1}\xi - x_{k,0}\xi^2| + |x_{k,0}\xi^2 - x_{k,2}| \leq (|\xi|+1)c_3 X_k^{-1}.$$

For $k \geq 2$, the recurrence formula of Lemma 2.6 (i) implies

(4.2)
$$x_{k+2,2} = a x_{k,0} x_{k+1,2} \pm x_{k-1,2}$$

and Lemma 2.6 (ii) gives
$$x_{k,0} x_{k+1,2} = x_{k,2} x_{k+1,0} \pm a x_{k-1,0} \pm 2 x_{k-1,1}.$$

Using (4.1), the latter relation leads to the estimate
$$\{x_{k,0} x_{k+1,2} \xi\} \leq X_k \{x_{k+1,0}\xi\} + a\{x_{k-1,0}\xi\} + 2\{x_{k-1,1}\xi\} \leq c_7 X_{k-1}^{-1}$$
for some constant $c_7 > 0$ (since $X_k X_{k+1}^{-1} \leq c_3^{2+\gamma} X_{k-1}^{-1}$ by virtue of Proposition 2.1 (i)). Combining this with (4.2), we deduce
$$|\delta_{k+2} - \delta_{k-1}| \leq a\{x_{k,0} x_{k+1,2} \xi\} \leq a c_7 X_{k-1}^{-1}.$$
Since the sequence $(X_k)_{k\geq 1}$ grows at least geometrically, this in turn implies that, for any pair of integers $j$ and $k$ which are congruent modulo 3 with $j \geq k \geq 1$, we have
$$|\delta_j - \delta_k| \leq c_8 X_k^{-1}$$
with some other constant $c_8 > 0$. Since
$$|\{x_{k,0}\xi^3\} - \delta_k| \leq |x_{k,0}\xi^3 - x_{k,2}\xi| \leq c_3|\xi|X_k^{-1}, \quad (k \geq 1),$$
we conclude that, for $i = 1, 2, 3$, the limit
$$\theta_i = \lim_{j\to\infty} \{x_{i+3j,0}\xi^3\} = \lim_{j\to\infty} \delta_{i+3j}$$
exists and that
$$|\theta_i - \{x_{k,0}\xi^3\}| \leq (c_8 + c_3|\xi|)X_k^{-1}$$
for $k \equiv i \bmod 3$. Since, for all sufficiently large $k$, Proposition 9.2 of [6] gives
$$\{x_{k,0}\xi^3\} \geq c_9 X_k^{-1/\gamma^3}$$
with a constant $c_9 > 0$, these numbers $\theta_i$ are non-zero. Thus the sequence $(\{x_{k,0}\xi^3\})_{k\geq 1}$ has (at most three) non-zero accumulation points and therefore is bounded below by some positive constant, say for $k \geq k_0$, to exclude the finitely many indices $k$ where $x_{k,0} = 0$. Applying Proposition 2.2 to the subsequence $(\mathbf{x}_k)_{k\geq k_0}$, we conclude that $\xi$ has the approximation property stated in Theorem 1.1.

Damien ROY
Département de Mathématiques et de Statistiques
Université d'Ottawa
585 King Edward
Ottawa, Ontario K1N 6N5, Canada
*E-mail:* droy@uottawa.ca
http://aix1.uottawa.ca/∼droy/